\documentclass[12pt,a4paper]{article}

\usepackage[T1]{fontenc}
\usepackage[utf8]{inputenc}
\usepackage{lmodern}
\usepackage[a4paper,margin=28mm]{geometry}
\usepackage{amsmath,amssymb,amsthm}
\usepackage{microtype}
\usepackage[hidelinks]{hyperref}

\newtheorem{theorem}{Theorem}[section]
\newtheorem{proposition}[theorem]{Proposition}
\newtheorem{lemma}[theorem]{Lemma}
\newtheorem{corollary}[theorem]{Corollary}
\theoremstyle{remark}
\newtheorem{remark}[theorem]{Remark}

\newcommand{\R}{\mathbb R}
\newcommand{\N}{\mathbb N}
\newcommand{\Nzero}{\mathbb N_0}
\newcommand{\Sym}{\operatorname{Sym}}
\newcommand{\tr}{\operatorname{tr}}
\newcommand{\spann}{\operatorname{span}}
\newcommand{\op}{\mathrm{op}}
\newcommand{\V}{\mathrm V}
\newcommand{\cT}{\mathcal T}
\newcommand{\norm}[1]{\left\lVert #1\right\rVert}
\newcommand{\errExt}{e_{\V}^{\mathrm{ext}}}
\newcommand{\errCube}{e_{\V}^{\square}}
\newcommand{\compExt}{n_{\V}^{\mathrm{ext}}}
\newcommand{\compCube}{n_{\V}^{\square}}
\newcommand{\errCinfty}{e_{F^\infty}^{\square}}
\newcommand{\compCinfty}{n_{F^\infty}^{\square}}

\title{Optimal Quadrature and the Curse of Dimensionality for a Veronese Class}
\author{Erich Novak \and Esra Sengelen Sevim}
\date{November 2026}

\hypersetup{
 pdftitle={Optimal Quadrature and the Curse of Dimensionality for a Veronese Class},
 pdfauthor={Erich Novak and Esra Sengelen Sevim},
 pdfsubject={Optimal exterior-node quadrature and information complexity},
 pdfcreator={LaTeX}
}

\begin{document}
\maketitle

\begin{center}
{\large Dedicated to Boris Kashin on the occasion of his 75th birthday} 
\end{center}

\begin{abstract}
We study integration on a $3^d$-dimensional
trigonometric space whose unit ball has all coordinate mixed partial
derivatives bounded by one. With evaluations anywhere in $\R^d$, we determine
the exact minimal worst-case error among all deterministic algorithms
using at most $n$ function values. An optimal quadrature formula
is obtained by truncating an exact product quadrature formula on
$\{1/2,1/2+\pi\}^d$. All weights are positive.
For each error tolerance $0<\varepsilon<1$, the information complexity
grows exponentially in $d$ with the common base $\Lambda\approx1.10546546$.
The correction of order $\sqrt d$ to the logarithm of the information
complexity and the corresponding error limit are determined by the
standard normal distribution function. For $b\ge1$ and
$n=\lceil b^d\rceil$, the error tends to one, one half, or zero
according as $b<\Lambda$, $b=\Lambda$, or $b>\Lambda$.
\end{abstract}

\section{Introduction and Main Results}

We study worst-case integration on a finite-dimensional trigonometric
class. The integral is taken over the unit cube, while function values
may be taken either in the cube or anywhere in $\R^d$.

For $d\in\N$, define
\[
 F_d^\infty=\left\{f\in C^\infty([0,1]^d;\R):
 \norm{D^\alpha f}_\infty\le1\ \text{for every }\alpha\in\Nzero^d\right\}.
\]
The supremum norm on $[0,1]^d$ is denoted by $\norm{\cdot}_\infty$, and
$D^0f=f$. The sets of positive and nonnegative integers are
$\N=\{1,2,\ldots\}$ and $\Nzero=\{0,1,2,\ldots\}$. All logarithms are
natural, empty sums equal zero, and empty products equal one.

Wo\'zniakowski \cite{Wozniakowski2003} posed the tractability question
for $F_d^\infty$. Wojtaszczyk \cite{Wojtaszczyk2003} proved that the
problem is not strongly polynomially tractable. Sengelen Sevim
\cite{Sevim2026} proved the curse of dimensionality, with asymptotic
exponential base approximately $1.10546546$ in the lower bound for the
information complexity. Matrix methods for quadrature in Hilbert
spaces were developed in \cite{HKKV2021}.

Here we show that the spectral lower bound from \cite{Sevim2026}
is attained on the Veronese class when exterior nodes are allowed.
We determine optimal quadrature formulas for every $n$ and derive
Gaussian asymptotics for their error and information complexity.

For $x\in\R$, define
\[
 z(x)=\begin{pmatrix}\cos((x-1/2)/2)\\ \sin((x-1/2)/2)\end{pmatrix},
 \qquad
 z_d(x_1,\ldots,x_d)=\bigotimes_{j=1}^d z(x_j).
\]
For $N\in\N$, let $\Sym_N$ denote the real symmetric $N\times N$ matrices,
and let $\norm{M}_{\op}$ be the operator norm induced by the Euclidean norm.
The Veronese class is
\[
 B_{\V,d}=\left\{f_M:\ M\in\Sym_{2^d},\ \norm{M}_{\op}\le1\right\},
 \qquad
 f_M(x)=z_d(x)^{\mathsf T}Mz_d(x).
\]
The linear span has dimension $3^d$. Define
\[
 S_d(f)=\int_{[0,1]^d}f(x)\,dx,
 \qquad D_d^{\mathrm{ext}}=\R^d,
 \qquad D_d^{\square}=[0,1]^d.
\]
For $\sigma\in\{\mathrm{ext},\square\}$ and $n\in\Nzero$,
$e_{\V}^{\sigma}(n,d)$ denotes the minimal deterministic worst-case error
among all algorithms based on at most $n$ function values in $D_d^\sigma$.
Adaptive nodes, input-dependent stopping, and nonlinear return maps are
included. Thus
\[
 e_{\V}^{\sigma}(n,d)=
 \inf_A\sup_{f\in B_{\V,d}}|S_d(f)-A(f)|.
\]
By Proposition~\ref{prop:embedding}, the initial errors satisfy
\[
 e_{\V}^{\mathrm{ext}}(0,d)=e_{\V}^{\square}(0,d)=1.
\]
For $0<\varepsilon<1$, define
\[
 n_{\V}^{\sigma}(\varepsilon,d)=
 \min\{n\in\Nzero:e_{\V}^{\sigma}(n,d)\le\varepsilon\}.
\]
The absolute and normalized error criteria coincide. Exponential
growth of $n_{\V}^{\sigma}(\varepsilon,d)$ in $d$ for a fixed
$\varepsilon\in(0,1)$ is called the curse of dimensionality.

Set
\[
 q=\frac12+\sin\frac12,
 \qquad p=\frac12-\sin\frac12,
 \qquad H(p)=-p\log p-q\log q,
 \qquad \Lambda=\exp(H(p)),
\]
and
\[
 V=pq\bigl[\log(q/p)\bigr]^2.
\]
Numerically,
\[
\begin{aligned}
 q&\approx0.97942554, & p&\approx0.02057446,\\
 H(p)&\approx0.10026647, & \Lambda&\approx1.10546546,
 & \sqrt V&\approx0.54835925.
\end{aligned}
\]
Let $\Phi$ be the standard normal distribution function, and let
$\Phi^{-1}$ denote the inverse function. For $d\in\N$, arrange the numbers
$q^{d-k}p^k$, with multiplicity $\binom dk$ for $k=0,\ldots,d$, in
decreasing order as
\[
 \lambda_1^{(d)}\ge\cdots\ge\lambda_{2^d}^{(d)}>0.
\]
For $n\in\Nzero$, define
\[
 \tau_{n,d}=
 \begin{cases}
  \displaystyle\sum_{j=n+1}^{2^d}\lambda_j^{(d)},&0\le n<2^d,\\[2mm]
  0,&n\ge2^d.
 \end{cases}
\]

\begin{theorem}\label{thm:main}
Let $d\in\N$ and $n\in\Nzero$. Then
\[
 \errCube(n,d)\ge\errExt(n,d)=\tau_{n,d}.
\]
For $n\le2^d$, $\errExt(n,d)$ is attained by a nonadaptive quadrature formula
with $n$ nodes from $\{1/2,1/2+\pi\}^d$. The weights are positive and equal
to the $n$ largest eigenvalues. For $n\ge2^d$, the complete grid attains
zero error with $2^d$ evaluations.
\end{theorem}

\begin{theorem}\label{thm:gaussian}
Let $s\in\R$, and let $(n_d)_{d\ge1}\subset\N$ satisfy
\[
 \frac{\log n_d-dH(p)}{\sqrt d}\longrightarrow s.
\]
Then
\begin{equation}\label{eq:gaussian-error}
 \lim_{d\to\infty}\errExt(n_d,d)
 =1-\Phi\!\left(\frac{s}{\sqrt V}\right).
\end{equation}
For every fixed $0<\varepsilon<1$,
\begin{equation}\label{eq:gaussian-complexity}
 \lim_{d\to\infty}
 \frac{\log\compExt(\varepsilon,d)-dH(p)}{\sqrt d}
 =\sqrt V\,\Phi^{-1}(1-\varepsilon).
\end{equation}
Consequently,
\begin{equation}\label{eq:complexity-main}
 \lim_{d\to\infty}\compExt(\varepsilon,d)^{1/d}=\Lambda.
\end{equation}
\end{theorem}

\begin{remark}[Monte Carlo]
For $n\in\N$ and independent uniform points $\Xi_1,\ldots,\Xi_n$ in $[0,1]^d$, put
$Q_{n,d}^{\mathrm{MC}}(f)=n^{-1}\sum_{i=1}^nf(\Xi_i)$.
Independence and $|f|\le1$ give
$\mathbb E|S_d(f)-Q_{n,d}^{\mathrm{MC}}(f)|^2
=\operatorname{Var}(f(\Xi_1))/n\le1/n$ for $f\in F_d^\infty$.
Thus
\begin{equation}\label{eq:MC}
 \sup_{f\in F_d^\infty}
 \bigl(\mathbb E|S_d(f)-Q_{n,d}^{\mathrm{MC}}(f)|^2\bigr)^{1/2}
 \le n^{-1/2}.
\end{equation}
The same bound applies to $B_{\V,d}$.
\end{remark}

\begin{remark}[Optimal Quadrature Formulas by Deletion]
The optimal exterior formulas are obtained from one exact product
formula by retaining its largest weights. A related deletion property
holds for the periodic trigonometric problem
\cite{Novak1999,Vybiral2020}; see Section~\ref{sec:three-problems}.
\end{remark}

\begin{remark}[Approximation, Integration, and Widths]
For $F_d^\infty$, let $e^{\mathrm{app}}(n,d)$ denote the minimal
deterministic worst-case error for uniform approximation by algorithms based
on at most $n$ continuous linear functionals, selected adaptively, and an
arbitrary return map. Novak and Wo\'zniakowski
\cite[Theorem~1]{NovakWozniakowski2009} proved
\[
 e^{\mathrm{app}}(n,d)=1
 \qquad\bigl(0\le n<2^{\lfloor d/2\rfloor}\bigr).
\]
The corresponding integration question is resolved for $F_d^\infty$ by the
exponential lower bound in \cite{Sevim2026}.

For widths of smooth function classes and related numbers, see
Kashin \cite{Kashin1977} and Pietsch \cite{Pietsch1987}; for multivariate
approximation and optimal sampling, see Temlyakov \cite{Temlyakov2018}
and the survey of Krieg and Ullrich \cite{KriegUllrich2026}.

Novak \cite[Proposition~3(c)]{Novak1986} constructed an infinite-dimensional
linear subspace $E\subset C([0,1])$ such that, for every continuous
linear functional $\mu$ on $C([0,1])$, there are $x_\mu\in[0,1]$
and $a_\mu\in\R$ with
\[
 \mu(f)=a_\mu f(x_\mu)\qquad(f\in E).
\]
Thus every such functional can be computed exactly on $E$ from one
function value. 
On the whole space $E$, the
Kolmogorov and Gelfand widths in the uniform norm are infinite for every
finite index, and uniform approximation from finitely many function
values has infinite worst-case error. Hence integration and uniform
approximation can have a very different information complexity; see
\cite{NovakWozniakowski2010}.
\end{remark}

\section{Matrix Representation and Optimal Quadrature}

For $t\in\R$, $x\in\R^d$, and $M\in\Sym_{2^d}$, define
\[
 K(t)=z(t)z(t)^{\mathsf T},\qquad
 K_d(x)=z_d(x)z_d(x)^{\mathsf T},\qquad
 (T_dM)(x)=\tr(MK_d(x)).
\]
Set
\[
 \cT_1=\spann\{1,\cos(x-\tfrac12),\sin(x-\tfrac12)\},
 \qquad \cT_d=\bigotimes_{j=1}^d\cT_1.
\]
For $f\in\cT_d$, define
\begin{equation}\label{eq:quotient-norm}
 \norm{f}_{\V,d}=
 \min_{\substack{M\in\Sym_{2^d}\\T_dM=f}}\norm{M}_{\op}.
\end{equation}
The class $B_{\V,d}$ is the unit ball of the quotient norm.

\subsection{Matrix Model and Derivative Bound}

Let $N\in\N$ and let $\mu_1,\ldots,\mu_N$ be the eigenvalues of
$B\in\Sym_N$. The trace norm, also called the nuclear norm, is defined by
$\norm{B}_*=\sum_{j=1}^N|\mu_j|$.

The spectral theorem yields
$\norm{B}_{\op}=\max_{1\le j\le N}|\mu_j|$ and, for every
$M\in\Sym_N$, the duality identities
\begin{equation}\label{eq:trace-pairing}
 |\tr(MB)|\le\norm{M}_{\op}\norm{B}_*,\qquad
 \max_{\substack{M\in\Sym_N\\\norm{M}_{\op}\le1}}|\tr(MB)|=\norm{B}_*.
\end{equation}
For an orthonormal eigenbasis $v_1,\ldots,v_N$ of $B$, we obtain
$\tr(MB)=\sum_{j=1}^N\mu_jv_j^{\mathsf T}Mv_j$ and
$|v_j^{\mathsf T}Mv_j|\le\norm{M}_{\op}$. Equality in
\eqref{eq:trace-pairing} follows for the symmetric matrix determined by
$Mv_j=v_j$ when $\mu_j\ge0$ and $Mv_j=-v_j$ when $\mu_j<0$.

For $N_1,N_2\in\N$, $A\in\Sym_{N_1}$, and $B\in\Sym_{N_2}$, every eigenvalue of
$A\otimes B$ is a product of an eigenvalue of $A$ and an eigenvalue of $B$.
Consequently,
\begin{equation}\label{eq:tensor-norms}
 \norm{A\otimes B}_{\op}=\norm{A}_{\op}\norm{B}_{\op},\qquad
 \norm{A\otimes B}_*=\norm{A}_*\norm{B}_*.
\end{equation}
The double-angle identities give
\begin{equation}\label{eq:K1}
 K(x)=\frac12
 \begin{pmatrix}
  1+\cos(x-\tfrac12)&\sin(x-\tfrac12)\\
  \sin(x-\tfrac12)&1-\cos(x-\tfrac12)
 \end{pmatrix}.
\end{equation}
Since $\norm{z(x)}_2=1$, the matrix $K(x)$ is a rank-one orthogonal
projection. The entries span $\cT_1$. For
$f(x)=c_0+c_1\cos(x-\tfrac12)+c_2\sin(x-\tfrac12)$, the identity
$f(x)=\tr(A_fK(x))$ determines a unique matrix
\begin{equation}\label{eq:A1}
 A_f=\begin{pmatrix}c_0+c_1&c_2\\c_2&c_0-c_1\end{pmatrix}.
\end{equation}
The eigenvalues of $A_f$ are $c_0\pm\sqrt{c_1^2+c_2^2}$. Hence
\begin{equation}\label{eq:V1}
 \norm{f}_{\V,1}=\norm{A_f}_{\op}
 =|c_0|+\sqrt{c_1^2+c_2^2}.
\end{equation}
The name Veronese refers to the quadratic map
$(u,v)\mapsto(\mathsf X,\mathsf Y,\mathsf Z)=(u^2,uv,v^2)$.
For the entries $u,v$ of $z(x)$, the image coordinates satisfy
$\mathsf X+\mathsf Z=1$ and
$\mathsf Y^2=\mathsf X\mathsf Z$.
Thus the quadratic image of $z$ lies on a conic in the plane
$\mathsf X+\mathsf Z=1$.
The matrix estimates depend on the rank-one property of $K(x)$.

The tensor definition gives

\[
 K_d(x)=\bigotimes_{j=1}^dK(x_j).
\]
Thus $K_d(x)$ is a rank-one orthogonal projection on $\R^{2^d}$. Every
entry of $K_d$ is in $\cT_d$. Conversely, tensor products of the matrices
in \eqref{eq:A1} yield every basis tensor of $\cT_d$, so
$T_d(\Sym_{2^d})=\cT_d$. The functions
$1,\cos(x-\tfrac12),\sin(x-\tfrac12)$ are linearly independent on
$[0,1]$: the value and the derivatives of orders one and two at $1/2$
determine all three coefficients. Tensor products therefore yield $3^d$
linearly independent functions on $[0,1]^d$.

The minimum in \eqref{eq:quotient-norm} is attained: a bounded minimizing
sequence in the closed affine space $\{M\in\Sym_{2^d}:T_dM=f\}$
has a convergent subsequence. Thus
$B_{\V,d}=\{f\in\cT_d:\norm{f}_{\V,d}\le1\}$.
For $d\ge2$, the representing matrix is not unique.
Indeed, $\dim\Sym_4=10>9=\dim\cT_2$, and tensoring a nonzero element
of $\ker T_2$ with $I_{2^{d-2}}$ gives a nonzero element of $\ker T_d$.
This explains the minimum in \eqref{eq:quotient-norm}.

Linear independence of the tensor basis on the cube implies that its
coefficients, and hence its extension to $\R^d$, are uniquely determined.
Thus exterior evaluation is well defined on the quotient space.

The univariate trace-norm identity on $[0,1]$, the cube-node spectral
lower bound, and the finite expression for $\tau_{n,d}$ are obtained in
\cite{Sevim2026}. Lemma~\ref{lem:kernel-derivative} yields the equality
on $\R$ and the tensor identity.

\begin{lemma}\label{lem:kernel-derivative}
Let $\ell\in\Nzero$ and $d\in\N$. For every $x\in\R$,

\[
 \norm{K^{(\ell)}(x)}_*=1.
\]
Consequently, for every $\alpha\in\Nzero^d$ and $x\in\R^d$,
\begin{equation}\label{eq:kernel-derivative-d}
 \norm{D^\alpha K_d(x)}_*=1.
\end{equation}

\end{lemma}

\begin{proof}
For $\ell=0$, the rank-one projection $K(x)$ has eigenvalues
one and zero. For $\ell\ge1$, differentiation of \eqref{eq:K1} yields
\[
 K^{(\ell)}(x)=\frac12
 \begin{pmatrix}
 \cos(x-\tfrac12+\tfrac{\ell\pi}{2})&\sin(x-\tfrac12+\tfrac{\ell\pi}{2})\\
 \sin(x-\tfrac12+\tfrac{\ell\pi}{2})&-\cos(x-\tfrac12+\tfrac{\ell\pi}{2})
 \end{pmatrix}.
\]
The characteristic polynomial equals $\lambda^2-\tfrac14$ because
$\cos^2t+\sin^2t=1$. The eigenvalues are $1/2$ and $-1/2$, and
$\norm{K^{(\ell)}(x)}_*=1$.

Each coordinate is present in one tensor factor, so

\[
 D^\alpha K_d(x)=\bigotimes_{j=1}^dK^{(\alpha_j)}(x_j).
\]
The tensor identity in \eqref{eq:tensor-norms} and the univariate trace-norm identity give
\eqref{eq:kernel-derivative-d}.
\end{proof}

For $f\in\cT_d$, define the coordinate derivative norm by

\[
 \norm{f}_{\infty,\mathrm{der}}
 =\sup_{\alpha\in\Nzero^d}\norm{D^\alpha f}_\infty.
\]

\begin{proposition}\label{prop:embedding}
Let $d\in\N$. For every $f\in\cT_d$,
\begin{equation}\label{eq:embedding}
 \norm{f}_{\infty,\mathrm{der}}\le\norm{f}_{\V,d}.
\end{equation}
Every $f\in B_{\V,d}$ satisfies $f|_{[0,1]^d}\in F_d^\infty$, and
$\errExt(0,d)=\errCube(0,d)=1$.
\end{proposition}

\begin{proof}
Let $\alpha\in\Nzero^d$ and $x\in[0,1]^d$. For every
$M\in\Sym_{2^d}$ satisfying $T_dM=f$, the trace estimate
\eqref{eq:trace-pairing} and Lemma~\ref{lem:kernel-derivative} imply

\[
 |D^\alpha f(x)|=|\tr(MD^\alpha K_d(x))|
 \le\norm{M}_{\op}\norm{D^\alpha K_d(x)}_*=\norm{M}_{\op}.
\]
The supremum over $x$ and $\alpha$, followed by the infimum over $M$,
yields \eqref{eq:embedding}. Since the derivative identity holds for
every $x\in\R^d$, we also obtain
\begin{equation}\label{eq:global-bound}
 \sup_{x\in\R^d}|D^\alpha f(x)|\le\norm{f}_{\V,d}
 \qquad(\alpha\in\Nzero^d).
\end{equation}
The identity matrix $I_{2^d}$ satisfies
$T_dI_{2^d}=1$, since $\tr K_d(x)=1$. Thus $\pm1\in B_{\V,d}$, and
$|S_d(f)|\le1$ for every $f\in B_{\V,d}$.
An algorithm with zero evaluations returns a constant $a\in\R$, and
$\max\{|1-a|,|-1-a|\}\ge1$. The constant approximation zero attains
error one.
\end{proof}

\subsection{Spectral Lower Bound}

The integration matrix and the spectrum are given in \cite{Sevim2026}.
Integration of \eqref{eq:K1} yields

\[
 G=\int_0^1K(x)\,dx=\begin{pmatrix}q&0\\0&p\end{pmatrix}.
\]
Coordinatewise integration yields

\[
 G_d=\int_{[0,1]^d}K_d(x)\,dx=G^{\otimes d},\qquad
 S_d(T_dM)=\tr(MG_d).
\]
The eigenvalues of $G_d$ are $\lambda_j^{(d)}$, $1\le j\le2^d$, and
$\tr(G_d)=(p+q)^d=1$.

For a deterministic algorithm $A$ and a class $B$, define
\[
 e(A,B)=\sup_{f\in B}|S_d(f)-A(f)|.
\]

\begin{proposition}\label{thm:spectral}
Let $d\in\N$ and $n\in\Nzero$. Then

\[
 \errCube(n,d)\ge\errExt(n,d)\ge\tau_{n,d}.
\]

\end{proposition}
\begin{proof}
Every cube-node algorithm is also an exterior-node algorithm, so the
first inequality is immediate. For $n\ge2^d$, the second follows from
$\tau_{n,d}=0$. Let $n<2^d$.
Run an arbitrary deterministic algorithm $A$ on the zero function.
Let $x^1,\ldots,x^m$, with $m\le n$, be the nodes it uses, and put
$a=A(0)$. Let $\Pi$ be the orthogonal projection onto
\[
 W=\spann\{z_d(x^1),\ldots,z_d(x^m)\},\qquad r=\dim W\le n,
\]
and set $P=I_{2^d}-\Pi$. Since $\norm P_{\op}\le1$, the functions
$\pm f_P$ belong to $B_{\V,d}$, where
\[
 f_P(x)=z_d(x)^{\mathsf T}Pz_d(x)=\norm{Pz_d(x)}_2^2.
\]
They vanish at every $x^i$. By induction over the evaluations, the
inputs $0,f_P,-f_P$ produce the same observations, the same stopping
decision and the same output $a$. Therefore
\[
 e(A,B_{\V,d})\ge
 \max\{|S_d(f_P)-a|,\ |-S_d(f_P)-a|\}\ge S_d(f_P).
\]
Let $v_1,\ldots,v_{2^d}$ be an orthonormal eigenbasis of
$G_d$, ordered by decreasing eigenvalues, and define
$b_j=\norm{\Pi v_j}_2^2$. Then $0\le b_j\le1$ and
$\sum_jb_j=\tr\Pi=r$. For $r\ge1$, we obtain
\[
 \begin{aligned}
 \tr(\Pi G_d)
 &=\sum_j\lambda_j^{(d)}b_j\\
 &=r\lambda_r^{(d)}+
   \sum_j\bigl(\lambda_j^{(d)}-\lambda_r^{(d)}\bigr)b_j\\
 &\le r\lambda_r^{(d)}+
   \sum_{j=1}^{r}\bigl(\lambda_j^{(d)}-\lambda_r^{(d)}\bigr)
 =\sum_{j=1}^{r}\lambda_j^{(d)}.
 \end{aligned}
\]
The inequality follows from $b_j\le1$ for $j\le r$ and $b_j\ge0$ for
$j>r$. For $r=0$, we have $\Pi=0$ and $\tr(\Pi G_d)=0$. Therefore
\[
 S_d(f_P)=\tr(PG_d)=1-\tr(\Pi G_d)
 \ge\sum_{j=r+1}^{2^d}\lambda_j^{(d)}\ge\tau_{n,d}.
\]
Taking the infimum over $A$ proves the exterior lower bound.
\end{proof}

\subsection{Optimal Exterior Quadrature}

For $u=(u_1,\ldots,u_d)\in\{0,1\}^d$, define
\begin{equation}\label{eq:externalgrid}
 \begin{aligned}
 |u|&=\sum_{j=1}^d u_j,\\
 x^u&=(1/2+\pi u_1,\ldots,1/2+\pi u_d),\\
 \lambda(u)&=q^{d-|u|}p^{|u|}.
 \end{aligned}
\end{equation}
Let $e_1=(1,0)^{\mathsf T}$ and $e_2=(0,1)^{\mathsf T}$.
For the nodes in \eqref{eq:externalgrid}, the identities $z(1/2)=e_1$
and $z(1/2+\pi)=e_2$ imply
$z_d(x^u)=\bigotimes_{j=1}^de_{1+u_j}$. Hence the vectors
$z_d(x^u)$ are the standard orthonormal eigenbasis of $G_d$.
With $P_u=K_d(x^u)$, the spectral decomposition becomes

\[
 G_d=\sum_{u\in\{0,1\}^d}\lambda(u)P_u.
\]
The trace pairing with a matrix satisfying $T_dM=f$ then yields a positive
quadrature formula on $\cT_d$.
\begin{equation}\label{eq:exterior-cubature}
 S_d(f)=\sum_{u\in\{0,1\}^d}\lambda(u)f(x^u),\qquad f\in\cT_d.
\end{equation}
For $d\in\N$ and $0\le n\le2^d$, let
$E_{n,d}\subseteq\{0,1\}^d$ contain $n$ elements satisfying
\[
 \lambda(u)\ge\lambda(v)\qquad
 (u\in E_{n,d},\ v\in\{0,1\}^d\setminus E_{n,d}).
\]
Define the quadrature formula
\begin{equation}\label{eq:optimal-exterior}
 Q_{n,d}^{\mathrm{ext}}(f)
 =\sum_{u\in E_{n,d}}\lambda(u)f(x^u),\qquad f\in\cT_d.
\end{equation}

\begin{proof}[Proof of Theorem~\ref{thm:main}]
Proposition~\ref{thm:spectral} provides
$\errCube(n,d)\ge\errExt(n,d)\ge\tau_{n,d}$.
For $f\in B_{\V,d}$, choose $M\in\Sym_{2^d}$ with $T_dM=f$ and
$\norm{M}_{\op}\le1$. Then, for $x\in\R^d$,
\[
 |f(x)|\le\norm{M}_{\op}\norm{z_d(x)}_2^2\le1.
\]
For $0\le n\le2^d$, equations~\eqref{eq:exterior-cubature}
and~\eqref{eq:optimal-exterior} therefore yield

\[
 |S_d(f)-Q_{n,d}^{\mathrm{ext}}(f)|
 \le\sum_{u\notin E_{n,d}}\lambda(u)=\tau_{n,d}.
\]
The function $f\equiv1\in B_{\V,d}$ attains equality, so

\[
 e(Q_{n,d}^{\mathrm{ext}},B_{\V,d})=\errExt(n,d)=\tau_{n,d}.
\]
The nodes are nonadaptive and the weights are positive. For $n=0$,
\eqref{eq:optimal-exterior} equals zero. For $n\ge2^d$, the full quadrature formula
in \eqref{eq:exterior-cubature} has zero error and uses at most $n$
evaluations.
\end{proof}

Define the cumulative multiplicities and masses by
\[
 L_{-1}(d)=B_{-1}(d)=0,\qquad
 L_r(d)=\sum_{k=0}^r\binom dk,\qquad
 B_r(d)=\sum_{k=0}^r\binom dk q^{d-k}p^k
\]
for $r=0,\ldots,d$.
For $n\in\{0,\ldots,2^d-1\}$, let $r\in\{0,\ldots,d\}$ satisfy
$L_{r-1}(d)\le n<L_r(d)$. Since $q>p$, the eigenvalue levels decrease
strictly with $k$. We obtain

\[
 \tau_{n,d}=(L_r(d)-n)q^{d-r}p^r
 +\sum_{k=r+1}^d\binom dk q^{d-k}p^k.
\]

\begin{corollary}\label{cor:finite-complexity}
Let $d\in\N$, $0<\varepsilon<1$, and
$r=\min\{j\in\{0,\ldots,d\}:B_j(d)\ge1-\varepsilon\}$. Then
\begin{equation}\label{eq:finite-complexity}
 \compExt(\varepsilon,d)
 =L_{r-1}(d)+
 \left\lceil
 \frac{1-\varepsilon-B_{r-1}(d)}{q^{d-r}p^r}
 \right\rceil.
\end{equation}
\end{corollary}

\begin{proof}
By Theorem~\ref{thm:main}, retain all levels below $r$, giving
$L_{r-1}(d)$ nodes and total weight $B_{r-1}(d)$, and then retain the
least integer number of weights $q^{d-r}p^r$ needed to reach total
weight $1-\varepsilon$. This gives \eqref{eq:finite-complexity}.
Since $B_{r-1}(d)<1-\varepsilon\le B_r(d)$, the ceiling lies between
$1$ and $\binom dr$.
\end{proof}

Since $q>p$, optimal node selection follows increasing $|u|$, with
arbitrary selection among equal weights. The sum of the weights is
$1-\tau_{n,d}$, and
\[
 Q_{n,d}^{\mathrm{ext}}(1)=1-\tau_{n,d}<1\qquad(0\le n<2^d).
\]

\begin{remark}[Trace-Norm Characterization]\label{thm:fixed-error}
For $Q_n(f)=\sum_{i=1}^n a_i f(x^i)$ with $x^i\in\R^d$ and
$a_i\in\R$, the definition of $B_{\V,d}$ and trace-norm duality give
\begin{equation}\label{eq:fixederror}
 e(Q_n,B_{\V,d})
 =\left\lVert G_d-\sum_{i=1}^n a_iK_d(x^i)\right\rVert_*.
\end{equation}
Indeed, the error is the maximum of
$|\tr(M[G_d-\sum_i a_iK_d(x^i)])|$ over $\norm M_{\op}\le1$.
This argument also covers noninjective $T_d$.

Since $B_{\V,d}$ is convex and symmetric and $S_d$ is a linear
functional, it suffices to consider nonadaptive linear quadrature
\cite{Bakhvalov1971,NovakWozniakowski2008}. Consequently,
\[
 \errCube(n,d)=
 \inf_{\substack{x^1,\ldots,x^n\in[0,1]^d\\a_1,\ldots,a_n\in\R}}
 \left\lVert G_d-\sum_{i=1}^n a_iK_d(x^i)\right\rVert_*.
\]
The same identity holds for $\errExt(n,d)$ with nodes in $\R^d$.
For $n=0$, the empty sum gives $\norm{G_d}_*=1$.
Formula~\eqref{eq:optimal-exterior} attains the exterior infimum for
$n\le2^d$; zero coefficients cover larger $n$.
\end{remark}

\section{Gaussian Limit and Information Complexity}

\begin{proof}[Proof of Theorem~\ref{thm:gaussian}]
Write
\[
 e_d(n)=\errExt(n,d).
\]
By Theorem~\ref{thm:main},
\[
 e_d(n)=1-\max_{\substack{I\subseteq\{0,1\}^d\\|I|\le n}}
 \sum_{u\in I}\lambda(u).
\]
Let $U$ satisfy $\mathbb P(U=u)=\lambda(u)$ and define
$Y_d=-\log\lambda(U)$. For every $n\in\N$ and $a>0$,
\begin{equation}\label{eq:counting-bound}
 \mathbb P(Y_d\ge\log n+a)-e^{-a}
 \le e_d(n)\le\mathbb P(Y_d\ge\log n).
\end{equation}
Fewer than $n$ points $u\in\{0,1\}^d$ satisfy $\lambda(u)>1/n$, which
yields the upper estimate.
For the lower estimate, every set $I$ with $|I|\le n$ satisfies
\[
 \mathbb P(U\in I,\ Y_d\ge\log n+a)
 =\sum_{\substack{u\in I\\-\log\lambda(u)\ge\log n+a}}\lambda(u)
 \le |I|\frac{e^{-a}}n\le e^{-a},
\]
so its omitted probability is at least
$\mathbb P(Y_d\ge\log n+a)-e^{-a}$.

Since $X_d=|U|$ is binomial with parameters $d$ and $p$,
\[
 Y_d=dH(p)+\log(q/p)(X_d-dp).
\]
The binomial central limit theorem gives
\[
 \frac{Y_d-dH(p)}{\sqrt{dV}}
 =\frac{X_d-dp}{\sqrt{dpq}}
 \Longrightarrow\mathcal N(0,1).
\]
Suppose that $(\log n_d-dH(p))/\sqrt d\to s$. Taking $a=d^{1/4}$ in
\eqref{eq:counting-bound}, both standardized thresholds converge to
$s/\sqrt V$, and $e^{-a}\to0$. Continuity of $\Phi$ yields
\eqref{eq:gaussian-error}.

Fix $0<\varepsilon<1$ and $\delta>0$. Define
\[
 s_\varepsilon=\sqrt V\,\Phi^{-1}(1-\varepsilon),\qquad
 m_d^\pm=\left\lceil
 \exp\bigl(dH(p)+(s_\varepsilon\pm\delta)\sqrt d\bigr)
 \right\rceil.
\]
Since $H(p)>0$, integer rounding changes the logarithms by $o(1)$.
Equation~\eqref{eq:gaussian-error} gives
\[
 e_d(m_d^-)>\varepsilon>e_d(m_d^+)
\]
for all sufficiently large $d$. Monotonicity yields
$m_d^-<\compExt(\varepsilon,d)\le m_d^+$.
Taking logarithms, subtracting $dH(p)$, dividing by $\sqrt d$,
and then letting $\delta\downarrow0$ gives \eqref{eq:gaussian-complexity}.
It follows that $\log\compExt(\varepsilon,d)/d\to H(p)$;
exponentiation gives \eqref{eq:complexity-main}.
\end{proof}

\begin{corollary}\label{cor:threshold}
Let $c>0$ and $b\ge1$ be independent of $d$. Then

\[
 \lim_{d\to\infty}\errExt\bigl(\lceil c b^d\rceil,d\bigr)
 =\begin{cases}
  1,&1\le b<\Lambda,\\[1mm]
  \tfrac12,&b=\Lambda,\\[1mm]
  0,&b>\Lambda.
 \end{cases}
\]

\end{corollary}
\begin{proof}
For $s\in\R$, put $\nu_d(s)=\lceil\exp(dH(p)+s\sqrt d)\rceil$.
Since $H(p)>0$, rounding changes the logarithm by $o(1)$, so
Theorem~\ref{thm:gaussian} gives
\begin{equation}\label{eq:gaussian-comparison}
 \lim_{d\to\infty}\errExt(\nu_d(s),d)=1-\Phi(s/\sqrt V).
\end{equation}
At $b=\Lambda$, Theorem~\ref{thm:gaussian} applies with $s=0$,
because $\log\lceil c\Lambda^d\rceil=dH(p)+O(1)$.

For $1\le b<\Lambda$, the relation
$\log\lceil cb^d\rceil=d\log b+O(1)$ implies
$\lceil cb^d\rceil\le\nu_d(s)$ eventually, for every fixed $s$.
Monotonicity and \eqref{eq:gaussian-comparison} give a lower limit
at least $1-\Phi(s/\sqrt V)$; letting $s\to-\infty$ gives one.
For $b>\Lambda$, the comparison is reversed, so the upper limit
is at most $1-\Phi(s/\sqrt V)$; letting $s\to+\infty$ gives zero.
Here we use $0\le\errExt(n,d)\le1$.
\end{proof}

\section{Cube Evaluations and Further Comparisons}

\subsection{Cube Evaluations}\label{sec:cube}

Let $\errCinfty(n,d)$ and $\compCinfty(\varepsilon,d)$ be defined in the same way,
with $F_d^\infty$ in place of $B_{\V,d}$ and nodes in $[0,1]^d$.
The embedding in Proposition~\ref{prop:embedding} and
$D_d^\square\subset D_d^{\mathrm{ext}}$ give
\[
 \errCinfty(n,d)\ge\errCube(n,d)\ge\errExt(n,d),\qquad
 \compCinfty(\varepsilon,d)\ge\compCube(\varepsilon,d)
 \ge\compExt(\varepsilon,d).
\]
Consequently, for every fixed $0<\varepsilon<1$,
\[
 \liminf_{d\to\infty}\compCube(\varepsilon,d)^{1/d}\ge\Lambda,
 \qquad
 \liminf_{d\to\infty}\compCinfty(\varepsilon,d)^{1/d}\ge\Lambda.
\]
For fixed $c>0$ and $1\le b<\Lambda$, both
$\errCube(\lceil cb^d\rceil,d)$ and
$\errCinfty(\lceil cb^d\rceil,d)$ tend to one.
At $b=\Lambda$, Theorem~\ref{thm:gaussian} and the comparison of
errors give
\[
 \liminf_{d\to\infty}\errCube(\lceil c\Lambda^d\rceil,d)\ge\tfrac12,
 \qquad
 \liminf_{d\to\infty}\errCinfty(\lceil c\Lambda^d\rceil,d)\ge\tfrac12.
\]
The exponential lower bound and the subcritical error limit for
$F_d^\infty$ were proved in \cite[Theorem~1.1]{Sevim2026}.

Define
\[
 \beta=\arccos(2\sin(1/2)),\qquad \xi_\pm=\tfrac12\pm\beta.
\]
Since $\cos(1/2)<2\sin(1/2)<1$, we have
$0<\beta<1/2$ and $0<\xi_-<\xi_+<1$. The sine terms in
\eqref{eq:K1} cancel at $\xi_-$ and $\xi_+$, and
$\cos\beta=2\sin(1/2)$. Therefore
\[
 G=\tfrac12K(\xi_-)+\tfrac12K(\xi_+).
\]
Tensorization gives an exact cube quadrature formula with $2^d$ nodes on
$\cT_d$. Since every eigenvalue of $G_d$ is positive,
Theorem~\ref{thm:main} implies
\[
 \errCube(n,d)=0\quad\Longleftrightarrow\quad
 \errExt(n,d)=0\quad\Longleftrightarrow\quad n\ge2^d.
\]

\subsection{Mixed and Directional Derivatives}
Coordinate mixed derivative bounds and Euclidean directional derivative
bounds define different classes. For $v=(v_1,\ldots,v_d)\in\R^d$, define
\[
 D_v=\sum_{j=1}^d v_j\frac{\partial}{\partial x_j}.
\]
Let $F_d^{\mathrm{dir}}$ be the class of all
$f\in C^\infty([0,1]^d;\R)$ such that
\[
 \sup_{x\in[0,1]^d}\sup_{v_1,\ldots,v_k\in\mathbb S^{d-1}}
 |D_{v_1}\cdots D_{v_k}f(x)|\le1
 \qquad(k\in\Nzero).
\]
Define
\[
 \mathcal P_d(x)=\prod_{j=1}^d x_j,\qquad
 v_d=d^{-1/2}(1,\ldots,1).
\]
Every coordinate mixed derivative of $\mathcal P_d$ has supremum norm at
most one, so $\mathcal P_d\in F_d^\infty$. In contrast,
\[
 \norm{D_{v_d}\mathcal P_d}_\infty=\sqrt d,
\]
and $\mathcal P_d\notin F_d^{\mathrm{dir}}$ for $d\ge2$. The comparison
agrees with the weak tractability results for directional derivative classes
\cite{HNU2014,HNUW2014b}.

\subsection{Veronese and Derivative Norms}\label{sec:norms}
For $1<\rho<\infty$, let $L_\rho=L_\rho([0,1])$ and let
$\norm{\cdot}_\rho$ denote its norm. Tensor product spaces with univariate norm
\[
 \left(\sum_{k=0}^{\infty}\norm{g^{(k)}}_\rho^\rho\right)^{1/\rho}
\]
are studied in \cite{NovakPillichshammer2025}.
The corresponding multivariate endpoint derivative norm is
$\sup_{\alpha\in\Nzero^d}\norm{D^\alpha f}_\infty$,
whose unit ball is $F_d^\infty$.

\begin{proposition}\label{prop:crossnorm}
Let $d\in\N$ and let $f_1,\ldots,f_d\in\cT_1$. Then
\[
 \norm{f_1\otimes\cdots\otimes f_d}_{\V,d}
 =\prod_{j=1}^d\norm{f_j}_{\V,1}.
\]
\end{proposition}

\begin{proof}
The tensor product of the matrices in \eqref{eq:A1}, together with
\eqref{eq:tensor-norms}, gives
\[
 \norm{f_1\otimes\cdots\otimes f_d}_{\V,d}
 \le\prod_{j=1}^d\norm{f_j}_{\V,1}.
\]
Write
\[
 f_j(x)=c_{0,j}+r_j\cos(x-\tfrac12-\theta_j).
\]
Equation~\eqref{eq:V1} gives
$\sup_{x\in\R}|f_j(x)|=\norm{f_j}_{\V,1}$.
The variables vary independently, and \eqref{eq:global-bound} gives the
reverse estimate.
\end{proof}

We show that the Veronese norm 
$\norm{f}_{\V,d}$ can be much larger than the norm $\norm{f}_{\infty,\mathrm{der}}$.
Define
\[
 C_d=\sup_{0\ne f\in\cT_d}
 \frac{\norm{f}_{\V,d}}{\norm{f}_{\infty,\mathrm{der}}}.
\]

\begin{proposition}
Let $d\in\N$. Then
\begin{equation}\label{eq:norm-gap}
 \left(2+\cos\tfrac12\right)^d\le C_d\le4^d.
\end{equation}
For $d=1$,
\begin{equation}\label{eq:C1}
 C_1=2+\cos\tfrac12.
\end{equation}
\end{proposition}

\begin{proof}
Put
\[
 f_*(x)=\cos(x-\tfrac12)-(1+\cos(1/2)).
\]
Then
\[
 \norm{f_*}_{\infty,\mathrm{der}}=1,
 \qquad \norm{f_*}_{\V,1}=2+\cos\tfrac12.
\]
Since $\norm{f_*^{\otimes d}}_{\infty,\mathrm{der}}=1$,
Proposition~\ref{prop:crossnorm} gives the lower estimate in
\eqref{eq:norm-gap}.

For the upper estimate, expand
\[
 f(x)=\sum_{\nu\in\{0,1,2\}^d}c_\nu
 \prod_{j=1}^d\chi_{\nu_j}(x_j-\tfrac12),
 \qquad \chi_0=1,\quad\chi_1=\cos,\quad\chi_2=\sin.
\]
Each basis tensor has $\V$-norm one. The coefficient $c_\nu$ is obtained
at $(1/2,\ldots,1/2)$ by applying in coordinate $j$ the operator
$\operatorname{Id}+\partial_j^2$, $-\partial_j^2$, or $\partial_j$ according
as $\nu_j=0,1,2$. Hence
\[
 \norm{f}_{\V,d}
 \le\sum_{\nu\in\{0,1,2\}^d}|c_\nu|
 \le4^d\norm{f}_{\infty,\mathrm{der}}.
\]
Let $d=1$ and normalize $\norm{f}_{\infty,\mathrm{der}}=1$. Write
\[
 f(x)=c_0+c_1\cos(x-\tfrac12)+c_2\sin(x-\tfrac12),
 \qquad a=|c_1|,\quad b=|c_2|.
\]
The values of $f'$ and $f''$ at $1/2$ give
$a,b\le1$. Put $c=\cos(1/2)$ and $t=\sin(1/2)$. The endpoint values give
\[
 |c_0+cc_1|+tb
 =\max\{|f(0)|,|f(1)|\}\le1.
\]
Equation~\eqref{eq:V1} therefore yields
\[
 \norm{f}_{\V,1}
 \le1+ca-tb+\sqrt{a^2+b^2}.
\]
The function on the right, without the initial one, is convex on
$[0,1]^2$, so the maximum is attained at a vertex. The four values are
$0$, $c+1$, $1-t$, and $c-t+\sqrt2$.
Since $t\ge\tfrac12-(\tfrac12)^3/6=23/48>\sqrt2-1$, the maximum is $c+1$.
Thus $\norm{f}_{\V,1}\le2+c$. The function $f_*$ gives equality in
\eqref{eq:C1}.
\end{proof}

\subsection{On optimal quadrature formulas}  \label{sec:three-problems}

For a nonempty class $\mathcal B\subset C([0,1]^d)$, write
\[
 e_n(\mathcal B)=
 \inf_{\substack{x^1,\ldots,x^n\in[0,1]^d\\a_1,\ldots,a_n\in\R}}
 \sup_{f\in\mathcal B}
 \left|S_d(f)-\sum_{i=1}^na_if(x^i)\right|.
\]
We know only a few cases where optimal quadrature formulas are known, 
in particular in the multivariate case, $d>1$. 

\medskip

1.  Let $B_d^{\mathrm{per}}$ be the unit ball of the $d$-fold Hilbert
tensor product of the three-dimensional space in which
$1,\cos(2\pi x),\sin(2\pi x)$ are declared orthonormal.
Its minimal linear error satisfies
\[
 e_n(B_d^{\mathrm{per}})^2=\max\{1-n2^{-d},0\}.
\]
Novak \cite{Novak1999} proved the lower bound for positive weights;
Vyb\'iral \cite[Section~5]{Vybiral2020} extended it to arbitrary real
weights. For $n\le2^d$, equality is attained by retaining any $n$
nodes of $\{0,1/2\}^d$, each with weight $2^{-d}$.

Thus both this periodic problem and the exterior Veronese problem admit
an exact $2^d$-point product formula whose summands can be deleted
without changing the remaining nodes or weights. In the periodic case
any $n$ summands may be retained; in the Veronese case the $n$ largest
weights must be retained, with arbitrary choices among ties.
This deletion property does not hold if we only allow nodes 
inside the cube. 

The Korobov-type class in \cite{SloanWozniakowski1997} is defined by a
weighted supremum bound on Fourier coefficients and has a different
norm. The Veronese ball is instead the image of an operator-norm
matrix ball; these norms define distinct integration problems.

\medskip

2.  For the H\"older class
\[
 \mathcal L_{\alpha,d}
 =\{f:[0,1]^d\to\R:
 |f(x)-f(y)|\le\norm{x-y}_\infty^\alpha
 \ \text{for all }x,y\in[0,1]^d\},
\]
where $0<\alpha\le1$, the quantization identity
\cite{Sukharev1979} gives, for $n\ge1$,
\[
 e_n(\mathcal L_{\alpha,d})
 =\inf_{\substack{X\subset[0,1]^d\\1\le|X|\le n}}
 \int_{[0,1]^d}\operatorname{dist}_\infty(x,X)^\alpha\,dx.
\]
If $n=m^d$ with $m\in\N$, the centers of the $m^d$ congruent subcubes
and equal weights $m^{-d}$ are optimal, and
\[
 e_{m^d}(\mathcal L_{\alpha,d})
 =\frac{d}{d+\alpha}(2m)^{-\alpha}.
\]

\medskip

3. 
Optimal univariate quadrature on Sobolev classes with derivative bounds in
$L_\rho$ can be characterized by minimal-norm monosplines. 
This was a very active research direction, see the excellent survey 
\cite{Zhensykbaev1981}.

\end{document}